\documentclass[11pt]{article}
\usepackage{amsmath,amssymb}
\newcommand{\doublespace}
   {\addtolength{\baselineskip}{0.15\baselineskip}}

\newcommand{\qed}{\hfill\rule{2mm}{2mm}\\}
\newtheorem{pdef}{Definition}[section] %
\newtheorem{thm}[pdef]{Theorem}        

\newcounter{equationnumber}
\renewcommand{\theequation}{\thesection.\arabic{equation}}
\def\mathletters{
    \addtocounter{equation}{1}
    \edef\@currentlabel{\theequation}
    \setcounter{equationnumber}{\value{equation}}
    \setcounter{equation}{0}
    \edef\theequation{\@currentlabel\noexpand\alph{equation}}
    }

\title{Characterization of the matrix whose norm is determined by
its action on decreasing sequences\\ (The exceptional cases)}
\author{Chang-Pao Chen$^\ast$, Chun-Yen Shen, and Kuo-Zhong Wang}

\date{ }
\begin{document}
\maketitle \doublespace
\pagestyle{myheadings} \thispagestyle{plain} \markboth{   }{ }
\begin{abstract}
Let $A=(a_{j,k})_{j,k \ge 1}$ be a non-negative matrix. In this
paper, we characterize those $A$ for which $\|A\|_{\ell_p,\ell_q}$
are determined by their actions on non-negative decreasing
sequences, where one of $p$ and $q$ is 1 or $\infty$. The
conditions forcing on $A$ are sufficient and they are also
necessary for non-negative finite matrices.
\end{abstract}
\footnotetext[1]{This work is supported by the National Science
Council, Taipei, ROC, under Grant NSC 94-2115-M-007-008}

\section{Introduction}

 For $x=\{x_k\}_{k=1}^\infty$, we write $x\ge 0$ if $x_k\ge 0$ for all $k$. We also write
$x\downarrow$ for the case that $\{x_k\}_{k=1}^\infty$ is
decreasing, that is, $x_k\ge x_{k+1}$ for all $k\ge 1$. For a
matrix $A=(a_{j,k})_{j,k\ge 1}$,  let $\|A\|_{E,F}$ denote the
norm of $A$ when $Ax=y$ defines an operator from $x\in E$ to $y\in
F$, where $(E, \|\cdot\|_E)$ and $(F, \|\cdot\|_F)$ are two normed
sequence spaces. More precisely, $
\|A\|_{E,F}=\sup_{\|x\|_E=1}\|Ax\|_F$. Clearly, $\|A\|_{E,F}\ge
\|A\|_{E,F,\downarrow}$, where
$$
\|A\|_{E, F,\downarrow}:=\sup_{\|x\|_E=1,x\ge 0,
x\downarrow}\|Ax\|_F.$$ The study of $\|A\|_{E,F}$ has a long
history in the literature and it goes back to the works of Hardy,
Copson, and Hilbert (cf. [10]). In [10, Theorem 326], Hardy proved
that $\|A\|_{\ell_p,\ell_p}=p/(p-1)$ for $1<p<\infty$, where
$A=(a_{j,k})_{j,k\ge 1}$ is the Ces\`aro matrix, defined by
$$
 a_{j,k}=\left\{
\begin{array}{ll}
  1/j& \mbox{if}\quad k\le j,\\
  \;\;\;\;0\;          & \mbox{if}\quad k>j.
\end{array} \right.
$$
This result can be restated in the following form, called the
Hardy inequality:
$$
  \sum_{j=1}^\infty \biggl |\frac 1j\sum_{k=1}^j x_k\biggr |^p\le \biggl (\frac
  p{p-1}\biggr )^p\sum_{k=1}^\infty |x_k|^p\qquad(x=\{x_k\}_{k=1}^\infty\in \ell_p).$$
For general $A$, some of the related results can be found in [1],
[3], [4], [6], [9], [14], and the references cited there. We also
refer the readers to [5], [15], and [16] for the integral setting.
As for the exceptional cases $p=1$ or $\infty$, the readers can
invoke [8], [11], [18], and others.

The question of determining $\|A\|_{E,F,\downarrow}$ was raised by
Bennett (cf. [1, page 422] and [3, page 422]). In [3, Problem
7.23], Bennett asked the following upper bound problem for the
case $E=F=\ell_p$: When does the equality
$\|A\|_{E,F}=\|A\|_{E,F,\downarrow}$ hold? This problem has been
partially solved by [1, page 422], [6, Lemma 2.4], and [12,
Theorem 2]. Recently, in [7], the present authors gave a more
general setting, which includes these as special cases. They
characterized $A$ and proved that $E$ and $F$ can be $\ell_p$,
$d(w, p)$, or $\ell_p(w)$, where $d(w,p)$ is the Lorentz sequence
space associated with non-negative decreasing weights $w_n$ and
$\ell_p(w)$ consists of all sequences $x=\{x_k\}_{k=1}^\infty$
such that
$$
 \|x\|_{\ell_p(w)}:=\left(\sum_{k=1}^\infty
|x_k|^pw_k\right)^{1/p}<\infty.$$ However,  the case
$F=\ell_\infty$ is excluded in [7]. The main purpose of this paper
is to deal with this case. In fact, we shall give a
characterization of $A$ for the case that $E=\ell_p$ and
$F=\ell_q$, where one of $p$ and $q$ is equal to 1 or $\infty$.
The details are given in \S2-\S3.

\section{The cases $p=1$  or  $\infty$}

In this section, we investigate the upper bound equality
$\|A\|_{\ell_p,\ell_q}=\|A\|_{\ell_p,\ell_q, \downarrow}$ for the
cases $p=1$ or $\infty$. The first main result is for $p=1$.

\begin{thm}
Let $1 \le q \le \infty$ and $A=(a_{j,k})_{j,k \ge 1}$ with
$a_{j,k} \ge 0.$ Suppose that $\|A\|_{\ell_1,\ell_q}<\infty$. Then
$(2.1)\Longleftrightarrow (2.2)\Longrightarrow (2.3)$, where
$$
  \biggl (\sum_{j=1}^\infty a_{j,1}^q \biggr )^{1/q} = \sup_{k\ge 1}\biggl (\sum_{j=1}^\infty a_{j,k}^q \biggr )^{1/q},\leqno{(2.1)}$$
$$
 \sup_{\|x\|_{\ell_1}=1}\|Ax\|_{\ell_q}=\max_{\|x\|_{\ell_1}=1, x\ge 0, x\downarrow} \|Ax\|_{\ell_q},\leqno{(2.2)}$$
$$
\|A\|_{\ell_1,\ell_q}=\|A\|_{\ell_1,\ell_q,\downarrow}.\leqno{(2.3)}
$$
If in addition, $a_{j,k}=0$ for $k>k_0$, where $k_0$ is a given
positive integer, then $(2.1)-(2.3)$ are equivalent.
\end{thm}

\begin{pf}
 By [8, Theorem 10] and [11, Eq. (15)], we know that
$\displaystyle \|A\|_{\ell_1,\ell_q}= \sup_{k\ge 1}\biggl
(\sum_{j=1}^\infty a_{j,k}^q \biggr )^{1/q}<\infty$. Combining
this with $(2.1)$, we obtain
$$
 \sup_{\|x\|_{\ell_1}=1}\|Ax\|_{\ell_q}=\|A\|_{\ell_1,\ell_q}=\biggl (\sum_{j=1}^\infty a_{j,1}^q \biggr )^{1/q}
 =\|Ae_1\|_{\ell_q},$$
 where $e_1=(1,0,\dots)$ is decreasing. Hence, $(2.1)\Longrightarrow (2.2)$. Assume that
$(2.2)$ holds. Then for some $x\ge 0$, we have $x\downarrow, $
$\|x\|_{\ell_1}=1$, and $\|Ax\|_{\ell_q}=\|A\|_{\ell_1,\ell_q}$.
For such an $x$, it follows from [8, Theorem 10] and [11, Eq.
(15)] that
$$
\|Ax\|_{\ell_q}=\|A\|_{\ell_1,\ell_q}=\sup_{k\ge 1}\biggl
(\sum_{j=1}^\infty a_{j,k}^q \biggr)^{1/q}=\sup_{k\ge 1}
S_k=M,\leqno{(2.4)}
$$ where $S_{k}=\left(\sum_{j=1}^\infty
a_{j,k}^q \right)^{1/q}$ and $M=\sup_{k\ge 1}S_{k}$.
 For $1 \le q < \infty$,
the function $f(t)=t^q$ is convex on $[0,\infty)$. Hence, by the
fact that $x_1+\cdots+x_n+\cdots=1$, we get
\begin{eqnarray*}
(2.5)\quad  \|Ax\|_{\ell_q}^q&=&(a_{1,1}x_1+a_{1,2}x_2+\cdots)^q+\cdots+(a_{n,1}x_1+a_{n,2}x_2+\cdots)^q+\cdots \\
                &\le& x_1S_1^q+x_2S_2^q+\cdots+x_nS_n^q+\cdots \le M^q.
\end{eqnarray*}
Putting $(2.4)-(2.5)$ together yields
$x_1S_1^q+x_2S_2^q+\cdots=M^q$, and consequently, $
x_1(M^q-S_1^q)+x_2(M^q-S_2^q)+\cdots=0.$ We know that $x\ge 0$,
$x\downarrow$, and $M^q-S_k^q\ge 0$ for all $k$. Therefore,
$M=S_1$, that is, $(2.1)$ holds.  This establishes the equivalence
$(2.1)\Longleftrightarrow (2.2)$ for the case $1\le q<\infty$. For
$q=\infty$, replace $(2.5)$ by
 \begin{eqnarray*}
(2.6)\qquad\qquad  \|Ax\|_{\ell_\infty}&=&\sup_{j\ge 1} (a_{j,1}x_1+a_{j,2}x_2+\cdots)\\
                &\le& x_1S_1+x_2S_2+\cdots+x_nS_n+\cdots \le
                M,\hspace{ 1in}
\end{eqnarray*}
and modify the proof between $(2.5)$ and $(2.6)$. Then we shall
get the equivalence $(2.1)\Longleftrightarrow (2.2)$ for
$q=\infty$. Clearly, $(2.2)\Longrightarrow (2.3)$. It remains to
prove the last conclusion. Assume that $a_{j,k}=0$ for $j\ge 1$
and $k>k_0$. We shall prove
$$
  \sup_{\|x\|_{\ell_1}=1,x\ge 0,x\downarrow} \|Ax\|_{\ell_q}\le \|Ay\|_{\ell_q}\leqno{(2.7)}$$
for some $y$ with $y\ge 0, y\downarrow$, and $\|y\|_{\ell_1}=1$.
If so, then $(2.3)$ implies $(2.2)$ and we are done. We have
$\|A\tilde x\|_{\ell_q}\ge \|Ax\|_{\ell_q}$, where $\tilde
x=(\tilde x_1,\tilde x_2,\cdots)$, $\tilde
x_1=x_1+(x_{k_0+1}+x_{k_0+2}+\cdots)$, $\tilde x_k=x_k$ for $1<
k\le k_0$, and $\tilde x_k=0$ otherwise. Hence, this substitution
does not loose the value of the left-hand side of $(2.7)$. Without
loss of generality, the sequences $x$ and $y$ in $(2.7)$ will be
assumed to be of the form $\tilde \xi= (\xi_1,\cdots,\xi_{k_0},
0,\cdots)$. Set $\xi^*= (\xi_1,\cdots,\xi_{k_0})$ and
$\|\xi^*\|_{\ell_1}=\sum_{k=1}^{k_0} |\xi_k|.$ We know that the
set $\Omega=\{\xi^*: \tilde \xi\ge 0, \tilde \xi\downarrow, \mbox
{and}\, \|\tilde \xi\|_{\ell_1}=1\}$ is a non-empty compact subset
of $\Bbb R ^{k_0}$ and the mapping $\tilde A: \Omega\mapsto \Bbb
R$ is continuous, where $\tilde A\xi^*=\|A\tilde \xi\|_{\ell_q}$.
Hence, the sequence $y$ involved in $(2.7)$ exists. This completes
the proof of Theorem 2.1. \qed
\end{pf}

We know that $\displaystyle \|A\|_{\ell_1,\ell_q}= \sup_{k\ge
1}\biggl (\sum_{j=1}^\infty a_{j,k}^q \biggr )^{1/q}$, so the
condition $\|A\|_{\ell_1,\ell_q}<\infty$ in Theorem 2.1 can be
replaced by the statement that the quantity on the right side of
$(2.1)$ is finite. For a finite matrix, $a_{j,k}=0$ for
$\max(j,k)>k_0$, where $k_0$ exists. Moreover,
$\|A\|_{\ell_1,\ell_q}<\infty$. Hence, $(2.1)-(2.3)$ in Theorem
2.1 are equivalent for this case. In general, $(2.3)$ does not
imply $(2.1)$. A counterexample is given by the matrix
$$
  A=\left(
    \begin{array}{cccccc}
      1& 0 & 0 & 0 & 0 & \cdots \\
      0& 1 & 0 & 0 & 0 & \cdots \\
      0& 1/4 & 1 & 0 & 0 & \cdots \\
      0& 0& 1/4 & 1 &0 &\cdots\\
      0& 0 & 1/9 & 1/4 & 1&\cdots\\
      0& 0 & 0 & 1/9 & 1/4 & \cdots \\
      \vdots & \vdots & \vdots & \vdots & \vdots & \ddots
    \end{array} \right).$$
For $x\ge 0$ with $\|x\|_{\ell_1}=1$, we have
$$
  \|Ax\|_{\ell_1} = \sum_{k=1}^\infty \biggl (\sum_{j=1}^\infty a_{j,k}\biggr )x_k
           \le  \frac {\pi^2}6 \biggl (\sum_{k=1}^\infty x_k\biggr ) =\frac {\pi^2}6.
$$
This implies $\|A\|_{\ell_1,\ell_1}\le \pi^2 /6$. On the other
hand, the choice $x_n=({1 \over n},{1 \over n},\dots,{1 \over
n},0,\dots)$ gives $x_n\ge 0, x_n\downarrow,\|x_n\|_{\ell_1}=1,$
and
\begin{eqnarray*}
  \|Ax_n\|_{\ell_1} &=& \frac 1n\biggl (\sum_{j=1}^\infty a_{j,1}+\cdots +\sum_{j=1}^\infty a_{j,n}\biggr )\\
           &=& \frac 1n\biggl (1 +(1+\frac 14)+\cdots+(1+\frac 14+\frac 19+\cdots+\frac 1{n^2})\biggr )\\
           &\longrightarrow & 1+\frac 14+\frac 19+\cdots=\frac {\pi^2}6\qquad \mbox{ as }\,n\rightarrow \infty.
\end{eqnarray*}
This leads us to $\displaystyle
 \|A\|_{\ell_1,\ell_1}=\pi^2/6=\sup_{\|x\|_{\ell_1}=1,x\ge 0, x\downarrow} \|Ax\|_{\ell_1},$
which says that $(2.3)$ holds for $q=1$. However, we can easily
see that $(2.1)$ is false for $q=1$.

The next theorem deals with the case $p=\infty$.

\begin{thm}
Let $1\le q \le \infty$ and $A=(a_{j,k})_{j,k \ge 1}$ with
$a_{j,k} \ge 0.$ Then
$$
 \|A\|_{\ell_\infty,\ell_q}=\biggl (\sum_{j=1}^\infty \biggl ( \sum_{k=1}^\infty
a_{j,k}\biggr )^q\biggr
)^{1/q}=\|A\|_{\ell_\infty,\ell_q,\downarrow}.\leqno{(2.8)}$$
\end{thm}

\begin{pf} Consider $1 \le q <\infty$. For $x\ge 0$ with
$\|x\|_{\ell_\infty}=1$, we have
$$
\|Ax\|_{\ell_q}=\biggl (\sum_{j=1}^\infty \biggl (
\sum_{k=1}^\infty a_{j,k}x_k\biggr )^q\biggr )^{1/q} \le \biggl
(\sum_{j=1}^\infty \biggl ( \sum_{k=1}^\infty a_{j,k}\biggr
)^q\biggr )^{1/q},$$
 and the right-hand side of the above inequality is attained by $x=(1,1,\dots)$. Therefore,
$(2.8)$ holds for $1 \le q <\infty$. As for $q= \infty$,
$$
\|Ax\|_{\ell_\infty}=\sup_{j\ge 1}\biggl (\sum_{k=1}^\infty
a_{j,k}x_{k}\biggr ) \le \sup_{j\ge 1}\biggl (\sum_{k=1}^\infty
a_{j,k}\biggr ),$$ where $x\ge 0$ and $\|x\|_{\ell_\infty}=1$.
Moreover, the choice $x=(1,1,\dots)$ gives
$\|Ax\|_{\ell_\infty}=\sup_{j\ge 1}\biggl (\sum_{k=1}^\infty
a_{j,k}\biggr )$. Hence, $(2.8)$ holds for $q=\infty$ and the
proof is complete. \qed
\end{pf}

From $(2.8)$ and the proof of Theorem 2.2, we see that $
\|A\|_{\ell_\infty,\ell_q}<\infty$ if and only if $\biggl
(\sum_{j=1}^\infty \biggl ( \sum_{k=1}^\infty a_{j,k}\biggr
)^q\biggr )^{1/q}<\infty$. Moreover, under this condition, the
following equality also holds:
$$
 \sup_{\|x\|_{\ell_\infty}=1}
\|Ax\|_{\ell_q}=\max_{\|x\|_{\ell_\infty}=1, x\ge 0, x\downarrow}
\|Ax\|_{\ell_q}.\leqno{(2.9)}$$

\section{The cases $q=1$ or $\infty$}

In this section, we investigate the upper bound equality for the
cases $q=1$ or $\infty$. Since $p=1$ or $\infty$ have been
examined in Theorems 2.1-2.2, we exclude these two cases in the
following, that is, we only consider the case $1<p<\infty$.

\begin{thm}
Let $1<p<\infty $ and $A=(a_{j,k})_{j,k \ge 1}$ with $a_{j,k} \ge
0$. Suppose that $\|A\|_{\ell_p,\ell_1}<\infty$. Then
$(3.1)\Longleftrightarrow (3.2)\Longleftrightarrow (3.3)$, where
$$
 \sum_{j=1}^\infty a_{j,k} \mbox{ is decreasing in k},\leqno{(3.1)}$$
$$
 \sup_{\|x\|_{\ell_p}=1} \|Ax\|_{\ell_1}=\max_{\|x\|_{\ell_p}=1,x\ge 0,x\downarrow} \|Ax\|_{\ell_1},\leqno{(3.2)}$$
$$
 \|A\|_{\ell_p,\ell_1}=
 \|A\|_{\ell_p,\ell_1,\downarrow}.\leqno{(3.3)}
 $$
\end{thm}

\begin{pf}
By [8, page 699, Corollary 1], we know that
$$
\|A\|_{\ell_p,\ell_1}=\biggl (\sum_{k=1}^\infty \biggl
(\sum_{j=1}^\infty a_{j,k}\biggr )^{p^*}\biggr
)^{1/p^*}<\infty,\leqno{(3.4)}$$ where $1/p+1/p^*=1$. Set $S_k=
\sum_{j=1}^\infty a_{j,k}$. Then $(3.1)$ says that
$\{S_k\}_{k=1}^\infty$ is decreasing.  Let $x=(x_1,x_2,\dots)$,
where $x_k=\lambda S_{k}^{p^*-1}$ and $\lambda=\biggl
(\sum_{k=1}^{\infty} S_{k}^{p^*}\biggr )^{-1/p}$. Then $x \ge 0,
x\downarrow$, $\|x\|_{\ell_p}=1$, and
$$
  \|Ax\|_{\ell_1}=\biggl (\sum_{k=1}^\infty S_k^{p^*}\biggr )^{1/p^*}
 =\biggl (\sum_{k=1}^\infty \biggl (\sum_{j=1}^\infty
a_{j,k}\biggr )^{p^*}\biggr )^{1/p^*}=\|A\|_{\ell_p,\ell_1}.$$
Hence, $(3.1)\Longrightarrow (3.2)$. Clearly,
$(3.2)\Longrightarrow (3.3)$. We claim that $(3.3)\Longrightarrow
(3.2)\Longrightarrow (3.1)$. Assume that $(3.3)$ holds. By
$(3.4)$, $ \|A\|_{\ell_p,\ell_1}=\biggl (\sum_{k=1}^\infty
S_k^{p^*}\biggr )^{1/p^*}$, and so there exists some
$x^n=(x^n_1,x^n_2,x^n_3,\cdots)\in\ell_p$ such that $x^n\geq 0$,
$x^n\downarrow$, $\|x^n\|_{\ell_p}=1$, and
$\|Ax^n\|_{\ell_1}\longrightarrow \biggl (\sum_{k=1}^\infty
S_k^{p^*}\biggr )^{1/p^*}$  as $n\rightarrow \infty.$ We know that
$\{x^n_k: n\ge 1\}\subset [0,1]$ for each $k$. By the ``diagonal
process" (cf. [17, Theorem 7.23]), without loss of generality, we
can further assume that for each $k\ge 1$, $x^n_k$ converges to
some $\tilde x_k$ as $n\to\infty$. Set $\tilde x=(\tilde
x_1,\tilde x_2,\cdots)$. Then $\tilde x\ge 0$ and $\tilde
x\downarrow$. We shall claim that $\|\tilde x\|_{\ell_p}=1$ and
$\|A\tilde x\|_{\ell_1}= \|A\|_{\ell_p,\ell_1}$. If so, $(3.2)$
follows. For any $m\ge 1$, we have
$$
 \bigg(\sum_{k=1}^m \tilde x_k^p\bigg)^{1/p}
 =\lim_{n\to\infty}  \bigg(\sum_{k=1}^m\bigg(x^n_k\bigg)^p\bigg)^{1/p}\le
 \lim_{n\to \infty} \|x^n\|_{\ell_p}=1,$$
which implies $\|\tilde x\|_{\ell_p}\le 1$. We shall prove
$\|\tilde x\|_{\ell_p}\ge 1$ and $\|A\tilde x\|_{\ell_1}=
\|A\|_{\ell_p,\ell_1}$ simultaneously. By definitions, $ \|A\tilde
 x\|_{\ell_1}=\sum_{k=1}^\infty S_k\tilde x_k$ and $
\|Ax^n\|_{\ell_1}=\sum_{k=1}^\infty S_kx^n_k$. For $m\ge 1$, it
follows from  the H\"older inequality that
\begin{eqnarray*}
\bigg|\sum_{k=1}^m S_k\tilde x_k-\sum_{k=1}^\infty S_kx^n_k\biggr
| &\le& \bigg|\sum_{k=1}^m S_k(\tilde
x_k-x^n_k)\bigg|+\|x^n\|_{\ell_p}\biggl (\sum_{k=m+1}^\infty
S_k^{p^*}\biggr )^{1/p^*}\\
&=& \bigg|\sum_{k=1}^m S_k(\tilde x_k-x^n_k)\bigg|+\biggl
(\sum_{k=m+1}^\infty S_k^{p^*}\biggr )^{1/p^*}.
\end{eqnarray*}
This implies
\begin{eqnarray*}
\sum_{k=1}^m S_k\tilde x_k &\ge&
\|Ax^n\|_{\ell_1}-\bigg|\sum_{k=1}^m S_k(\tilde
x_k-x^n_k)\bigg|-\biggl
(\sum_{k=m+1}^\infty S_k^{p^*}\biggr )^{1/p^*}\\
&\longrightarrow&\biggl (\sum_{k=1}^\infty S_k^{p^*}\biggr
)^{1/p^*}-\biggl (\sum_{k=m+1}^\infty S_k^{p^*}\biggr
)^{1/p^*}\qquad \mbox{ as }\quad n\rightarrow \infty.
\end{eqnarray*}
Taking $m\to\infty$, we get $\|A\tilde x\|_{\ell_1}\ge \biggl
(\sum_{k=1}^\infty S_k^{p^*}\biggr )^{1/p^*}$. For the reverse
inequality, by the H\"older inequality and $\|\tilde
x\|_{\ell_p}\le 1$, we obtain
$$
\|A\tilde x\|_{\ell_1}=\sum_{k=1}^\infty S_k\tilde x_k \le
\|\tilde x\|_{\ell_p}\biggl(\sum_{k=1}^\infty S_k^{p^*}\biggr
)^{1/p^*}\le\biggl(\sum_{k=1}^\infty S_k^{p^*}\biggr
)^{1/p^*}.\leqno{(3.5)}$$ Putting these inequalities together, we
find that $\|A\tilde x\|_{\ell_1}= \biggl (\sum_{k=1}^\infty
S_k^{p^*}\biggr )^{1/p^*}=\|A\|_{\ell_p,\ell_1}$ and $\|\tilde
x\|_{\ell_p}=1.$ This finishes the proof of the implication:
$(3.3)\Longrightarrow (3.2)$. In fact, we get more. Since the
inequality signs in $(3.5)$ are equality signs. By the H\"older
inequality, we infer that $(\tilde x_1^p,\tilde x_2^p,\dots)$ and
$(S_1^{p^*},S_2^{p^*},\dots)$ are proportional. Since $\tilde
x_1^p\ge \tilde x_2^p\ge\cdots$, the sequence
$\{S_k^{p^*}\}_{k=1}^\infty$ is decreasing. This leads us to
$(3.1)$. We complete the proof. \qed
\end{pf}

We know that $\displaystyle \|A\|_{\ell_p,\ell_1}=\biggl
(\sum_{k=1}^\infty \biggl (\sum_{j=1}^\infty a_{j,k}\biggr
)^{p^*}\biggr )^{1/p^*}$. Hence, the condition
$\|A\|_{\ell_p,\ell_1}<\infty$ in Theorem 3.1 can be replaced by
$\biggl (\sum_{k=1}^\infty \biggl (\sum_{j=1}^\infty a_{j,k}\biggr
)^{p^*}\biggr )^{1/p^*}<\infty$. As Theorems 2.1-2.2 indicate,
Theorem 3.1 is false for the cases that $p=1$ or $\infty$.

In [7], the present authors indicate that the matrix $A$, defined
by $a_{2,2}=1$ and 0 otherwise, possesses the property:
$\|A\|_{\ell_p,\ell_\infty}>\|A\|_{\ell_p,\ell_\infty,\downarrow}$,
where $1\le p<\infty$. This phenomenon can be interpreted by
applying the following result to the case $\Lambda=\{2\}$.

\begin{thm}
Let $1<p< \infty, 1/p+1/p^*=1,$ and $A=(a_{j,k})_{j,k \ge 1}$ with
$a_{j,k} \ge 0$. Suppose that there exists a nonempty finite set
$\Lambda$ of positive integers with
$$
  \sup_{j\notin \Lambda} \biggl(\sum_{k=1}^\infty a_{j,k}^{p^*}\biggr )^{1/p^*}
  <\sup_{j\in \Lambda} \biggl(\sum_{k=1}^\infty a_{j,k}^{p^*}\biggr)^{1/p^*}<\infty.\leqno{(3.6)}$$
Then $(3.7)\Longleftrightarrow (3.8)\Longleftrightarrow (3.9)$,
where
\begin{itemize}
\item[$(3.7)$]\quad there exists some $l\in \Lambda$ such that
 $a_{l,1} \ge a_{l,2}\ge \cdots \ge a_{l,n}\ge \cdots$  and
 $\displaystyle
  \biggl (\sum_{k=1}^\infty a_{l,k}^{p^*}\biggr )^{1/p^*}
  =\sup_{j\ge 1} \biggl(\sum_{k=1}^\infty a_{j,k}^{p^*}\biggr)^{1/p^*},$
  \end{itemize}
$$
 \sup_{\|x\|_{\ell_p}=1}\|Ax\|_{\ell_\infty}=\max_{\|x\|_{\ell_p}=1,x\ge 0,x\downarrow} \|Ax\|_{\ell_\infty},\leqno{(3.8)}$$
$$
\|A\|_{\ell_p,\ell_\infty}=\|A\|_{\ell_p,\ell_\infty,\downarrow}.\leqno{(3.9)}$$
For the implication from $(3.7)$ to any of $(3.8)$ or $(3.9)$, the
condition that $\Lambda$ is finite is unnecessary.
\end{thm}

\begin{pf}
Putting  the Hellinger-Toeplitz theorem (see [2, page 29]), [8,
Theorem 10], and $(3.6)$ together, we obtain
$$
 \|A\|_{\ell_p,\ell_\infty}=\|A^t\|_{\ell_1,\ell_{p^*}}=\sup_{j\ge
1}\biggl(\sum_{k=1}^\infty
a_{j,k}^{p^*}\biggr)^{1/p^*}<\infty,\leqno{(3.10)}$$ where $A^t$
is the transpose of $A$. Assume that $(3.7)$ holds. Set
$x=(x_1,x_2,\dots)$, where $x_k=\lambda a_{l,k}^{p^*-1}$ and
$\lambda=\biggl (\sum_{k=1}^{\infty} a_{l,k}^{p^*}\biggr)^{-1/p}$.
Then $x \ge 0, x\downarrow$, $\|x\|_{\ell_p}=1$, and
$$
  \|Ax\|_{\ell_\infty}\ge \sum_{k=1}^\infty a_{l,k}x_k=\biggl (\sum_{k=1}^{\infty} a_{l,k}^{p^*}\biggr )^{1/p^*}. \leqno{(3.11)}$$
By $(3.7)$ and $(3.10)$, we get $
  \|Ax\|_{\ell_\infty}\ge  \|A\|_{\ell_p,\ell_\infty}$. This leads us to $(3.8)$. Clearly, $(3.8)\Longrightarrow (3.9)$.
  In the above argument,
  the assumption that $\Lambda$ is finite is unnecessary. We
  claim that $(3.9)\Longrightarrow (3.7)$. Assume that $(3.9)$
  holds. We know that $\Lambda$ is a finite set. Without loss of generality, we can assume that
$\biggl(\sum_{k=1}^\infty a_{r,k}^{p^*}\biggr)^{1/p^*}=\sup_{j\in
\Lambda} \biggl(\sum_{k=1}^\infty a_{j,k}^{p^*}\biggr)^{1/p^*}$
for all $r\in\Lambda$. Let $x\ge 0$, $x\downarrow$,
$\|x\|_{\ell_p}=1$, and $\|Ax\|_{\ell_\infty}>\gamma$, where
$\gamma=\sup_{j\notin \Lambda} \biggl(\sum_{k=1}^\infty
a_{j,k}^{p^*}\biggr)^{1/p^*}$. We have $\displaystyle
\|Ax\|_{\ell_\infty}=\sup_{j\ge 1} \biggl ( \sum_{k=1}^\infty
a_{j,k}x_k\biggr )$. For $j\ge 1$, the H\"older inequality implies
$$
  \sum_{k=1}^\infty a_{j,k}x_k \le \|x\|_{\ell_p}\biggl (\sum_{k=1}^\infty a_{j,k}^{p^*}\biggr )^{1/p^*}
  =\biggl (\sum_{k=1}^\infty a_{j,k}^{p^*}\biggr )^{1/p^*},$$
which gives $\displaystyle \sup_{j\notin \Lambda} \biggl (
\sum_{k=1}^\infty a_{j,k}x_k\biggr )\le  \sup_{j\notin \Lambda}
\biggl (\sum_{k=1}^\infty a_{j,k}^{p^*}\biggr )^{1/p^*}=\gamma.$
Thus, $\displaystyle \|Ax\|_{\ell_\infty}=\sum_{k=1}^\infty
a_{r,k}x_k$ for some $r\in\Lambda$. Since $\Lambda$ is a finite
set, we can find some $l\in\Lambda$ such that
$$
\sup_{\|x\|_{\ell_p}=1, x\ge 0, x\downarrow} \sum_{k=1}^\infty
a_{l,k}x_k=\|A\|_{\ell_p,\ell_\infty,\downarrow}.\leqno{(3.12)}$$
Putting $(3.6)$, $(3.9)$, $(3.10)$, and $(3.12)$ together yields
$$
\sup_{\|x\|_{\ell_p}=1, x\ge 0, x\downarrow} \sum_{k=1}^\infty
a_{l,k}x_k=\biggl(\sum_{k=1}^\infty
a_{l,k}^{p^*}\biggr)^{1/p^*},$$ which can be written in the form:
$\|\tilde A\|_{\ell_p,\ell_1}=\|\tilde
A\|_{\ell_p,\ell_1,\downarrow}$. Here $\tilde A=(\tilde
a_{j,k})_{j,k\ge 1}$ is  defined by $\tilde a_{l,k}=a_{l,k}$ and
$\tilde a_{j,k}=0$ for $j\neq l$. By Theorem 3.1, we get $(3.7)$.
The proof is complete.\qed
\end{pf}

From $(3.10)$, we see  that condition $(3.6)$ implies
$\|A\|_{\ell_p,\ell_\infty}<\infty$. It is clear that this
condition is automatically satisfied by any finite non-negative
matrix $A$. Applying Theorem 3.2 to this case, we find that
$(3.7)-(3.9)$ are equivalent for such kind of matrices. In
general, $(3.6)$ can not be taken off. The following matrix
provides us a counterexample:
$$
  A=(a_{j,k})_{j,k\ge 1}=\left(
    \begin{array}{cccccc}
      1/2& 1 & 1/3 & 1/4 & 1/5 & \cdots \\
      1& 0 & 0 & 0 & 0 & \cdots \\
      1& 1/2 & 0 & 0 & 0 & \cdots \\
      1& 1/2& 1/3 & 0 &0 &\cdots\\
      1& 1/2 & 1/3 & 1/4 & 0&\cdots\\
      1& 1/2 & 1/3 & 1/4 & 1/5 & \cdots \\
      \vdots & \vdots & \vdots & \vdots & \vdots & \ddots
    \end{array} \right).$$
Clearly, both of $(3.6)-(3.7)$ are not satisfied by any finite set
$\Lambda$. Let $x_k=(1/k)^{p^*-1}\biggl (\sum_{s=1}^n
(1/s)^{p*}\biggr )^{-1/p}$ for $1\le k\le n$ and 0 otherwise,
where $1<p<\infty$. Then
 $x\ge 0, x\downarrow,$ $\|x\|_{\ell^p}=1$, and
$\displaystyle
 \|Ax\|_{\ell_\infty}\ge \sum_{k=1}^n x_k/k=\biggl
(\sum_{k=1}^n (1/k)^{p^*}\biggr )^{1/p^*}$ for $n\ge 2$.
 This leads us to
$$
\|A\|_{\ell_p,\ell_\infty,\downarrow}\ge
\|Ax\|_{\ell_\infty}\ge\biggl (\sum_{k=1}^n (1/k)^{p^*}\biggr
)^{1/p^*}=\biggl (\sum_{k=1}^\infty a_{n+1,k}^{p^*}\biggr
)^{1/p^*}\qquad(n\ge 2).$$ Putting this with $(3.10)$ and letting
$n\to\infty$, we obtain $\|A\|_{\ell_p,\ell_\infty,\downarrow}\ge
\|A\|_{\ell_p,\ell_\infty}$. Hence, $(3.9)$ holds.

\end{document}